% The following is the AMSTeX-file of Flicker-Scheiderer-Sujatha:
% ``GROTHENDIECK'S THEOREM ON NON-ABELIAN $H^2$ AND LOCAL-GLOBAL PRINCIPLES''
% accepted for publication by the Journal of the AMS.

% fss.tex    March 6, 1998

\input amstex.tex
\loadmsam
\loadmsbm
\loadbold
\input amssym.tex
\overfullrule=0pt
\baselineskip=13pt plus 2pt

\documentstyle{amsppt}
\pageheight{46pc}
\pagewidth{33pc}
\magnification=1200
\parskip 2pt
\NoBlackBoxes

\def\n{\noindent}

\def\>{\,}

\def\Aut{\operatorname{Aut}}
\def\cd{\operatorname{cd}}
\def\ch{\operatorname{char}}
\def\comp{\mathbin{\scriptstyle\circ}} 
\def\Gal{\operatorname{Gal}}
\def\Int{\operatorname{Int}}
\def\int{\operatorname{int}}
\def\Inv{\operatorname{Inv}}
\def\ol{\overline}
\def\obs{\operatorname{obs}}
\def\Out{\operatorname{Out}}
\def\SAut{\operatorname{SAut}}
\def\SOut{\operatorname{SOut}}
\def\Spec{\operatorname{Spec}}
\def\Sper{\operatorname{Sper}}
\def\vcd{\operatorname{vcd}}

\def\rightmap#1{\mathrel{\mathop{\to}\limits^{#1}}}
\def\isomap{\mathrel{\buildrel\sim\over\to}} 
\def\dirlim#1{{\mathsurround=0pt\setbox0=\hbox{\rm lim}
  \mathchoice
  {\mathop{\vphantom{\big|}\smash{\hbox{\rlap{\copy0}\lower
  5.5pt\hbox to\wd0{\rightarrowfill}}}}_{#1}}
  {\mathop{\hbox{\rlap{\copy0}\lower5.5pt\hbox{\hbox to\wd0
  {\rightarrowfill}$\scriptstyle #1$}}}}{}{}}}

\def\lab#1 {\medbreak\n #1\enspace\ignorespaces }
\def\sqr#1#2{{\vcenter{\vbox{\hrule height.#2pt
     \hbox{\vrule width.#2pt height#1pt\kern#1pt\vrule
     width.#2pt}\hrule height.#2pt}}}}
\def\square{\mathchoice\sqr54\sqr54\sqr{3.5}3\sqr{2.5}3}
\def\Proof{\medbreak\n{\it Proof.}\enspace\ignorespaces}
\def\endproof{\unskip\nobreak\kern4pt\nobreak\hfill\discretionary
  {}{\hbox to\hsize{\hfill$\square$}}{\vbox{\hbox{$\square$}}}%
  \bigbreak} % (!)
\outer\def\beginsect#1\par{\vskip8pt plus.15\vsize\penalty-250
     \vskip0pt plus-.15\vsize\bigskip\vskip\parskip\message{#1}
     \leftline{\bf #1}\nobreak\bigskip\n}

\def\ab{\operatorname{ab}}
\def\red{\operatorname{red}}

\def\ad{\operatorname{ad}}
\def\ss{\operatorname{ss}}
\def\sc{\operatorname{sc}}
\def\olgss{{\ol G}^{\ss}}
\def\olgsc{{\ol G}^{\sc}}
\def\olzss{{\ol Z}^{\ss}}
\def\olzsc{{\ol Z}^{\sc}}
\def\gss{G^{\ss}}
\def\gsc{G^{\sc}}
\def\zss{Z^{\ss}}
\def\zsc{Z^{\sc}}
\rightheadtext{GROTHENDIECK'S THEOREM ON NON-ABELIAN $H^2$}
\topmatter
\title GROTHENDIECK'S THEOREM ON NON-ABELIAN$\,\, H^2$ AND LOCAL-GLOBAL 
PRINCIPLES
\endtitle
\author Yuval Z.~Flicker, Claus Scheiderer, R.~Sujatha\endauthor
\footnote"~"{1991 Mathematics Subject Classification: 14L30, 11R34, 12G05.}
\abstract A theorem of Grothendieck asserts that over a perfect field 
$k$ of cohomological dimension one, all non-abelian $H^2$-cohomology 
sets of algebraic groups are trivial. The purpose of this paper is to
establish a formally real generalization of this theorem. The
generalization --- to the context of perfect fields of {\it virtual\/}
cohomological dimension one --- takes the form of a local-global
principle for the $H^2$-sets with respect to the orderings of the
field. This principle asserts in particular that an element in
$H^2$ is neutral precisely when it is neutral in the real closure
with respect to every ordering in a dense subset of the real
spectrum of $k$. Our techniques provide a new proof of Grothendieck's
original theorem. An application to homogeneous spaces over $k$ is
also given.
\endabstract
\endtopmatter
\document

\subhead 0. Introduction\endsubhead

\n Let $k$ be a field with separable closure $k_s$ and Galois group
$\Gamma_k=\Gal(k_s/k)$. We assume that $k$ is perfect throughout
this introduction; then $k_s$ is an algebraic closure of $k$. Let 
$G$ be an algebraic group over $k$. (By this we always mean a smooth 
group scheme of finite type over $\Spec k$.) Recall that a homogeneous 
$G$-space is a nonempty smooth algebraic scheme
$X$ together with a (right) action by $G$ on $X$, both $X$ and the
action being defined over $k$, such that the action of $G(k_s)$ on
$X(k_s)$ is transitive. The $G$-space $X$ is called a principal
homogeneous $G$-space (or a $G$-torsor) if in addition this last
action is free, i.e.\ has trivial stabilizer subgroups. The
non-abelian Galois cohomology set $H^1(k,G)=H^1(\Gamma_k,G(k_s))$
is the set of isomorphism classes of principal homogeneous
$G$-spaces over $k$. A well known theorem of Steinberg
(``Conjecture~I'' of Serre; see [S2], III, \S 2.3, Thm.~1') asserts
that when $k$ has cohomological dimension at most one
($\cd(k)\le1$) and $G$ is connected and linear, the cohomology set
$H^1(k,G)$ is trivial. In other words, each principal homogeneous
$G$-space is trivial, which means, has a $k$-point.

There is a general definition of the second non-abelian cohomology
set in terms of gerbes (by Grothendieck, Dedecker, Giraud [G]; see 
also Deligne-Milne [DM], Breen [Br]). Springer [Sp] constructed a
non-abelian $H^2$-set in terms of group extensions; this set has an
equivalent description in terms of 2-cocycles. His approach was
recently taken up again by Borovoi [B1]. Here is a brief review of
this setup: Given an algebraic group $\ol G$ over $k_s$, let
$\SOut(\ol G/k)$ be the quotient of the group $\SAut(\ol G/k)$ of
$k$-semilinear automorphisms of $\ol G$ by the subgroup $\Int(\ol G)$ of
inner automorphisms of $\ol G$. A $k$-kernel in $\ol G$ is a
homomorphic section $\kappa$ of the natural homomorphism $\SOut(\ol
G/k)\to\Gamma_k$ which satisfies a certain continuity condition.
(The notion of continuity is delicate.)
% Both Springer and Borovoi failed here! Let's pray that we are
% doing it right now!
The set $H^2(k,\ol G,\kappa)$ classifies group extensions of $\ol
G(k_s)$ by $\Gamma_k$ 
compatible with the
kernel $\kappa$. It may be
empty. If it is not, it is a principal homogeneous set under the
abelian group $H^2(k,Z)$, where $Z$ is the center of $\ol G$,
considered as an algebraic group over $k$ by means of $\kappa$. The
elements of $H^2(k,\ol G,\kappa)$ which correspond to split
extensions are called neutral. The set of neutral elements is denoted
by $N^2(k,\ol G,\kappa)$, and the kernel $\kappa$ is said to
be trivial if there exists a neutral element. Theorem 3.5 of [Sp]
--- attributed to Grothendieck --- asserts that if $k$ is perfect
with $\cd(k)\le1$ and
$\ol G$ is any algebraic group over $k_s$, then any $k$-kernel
$\kappa$ in $\ol G$ is trivial, and $H^2 (k,\ol G,\kappa)$ consists
of a single element which is (therefore) neutral.

The interest in non-abelian $H^2$-cohomology stems at least in part
from its relation to homogeneous spaces. Let $G$ be an algebraic
group over $k$ and $X$ a homogeneous $G$-space, defined over $k$.
To $X$ one associates in a natural way a $k$-kernel $\kappa_X$ in
$\ol H$, where $\ol H$ is the stabilizer of a $k_s$-point of $X$,
and a class $\alpha_X$ in $H^2(k,\ol H,\kappa_X)$. The class $\alpha_X$
is neutral if and only if $X$ is dominated (over $k$) by some {\it
principal\/} homogeneous $G$-space. This fact can be regarded as
part of an ``exact sequence'' relating the $H^i$-sets of $G$ and
$\ol H$ to relative $H^i$-sets of $G$ mod~$\ol H$, $i=0,1,2$; see
[Sp], Prop.~1.27.

Grothendieck's theorem therefore implies Springer's theorem, which
says: Over a (perfect) field $k$ with $\cd(k)=1$, each homogeneous
space under an algebraic group $G$ is dominated by a principal
homogeneous space under $G$. Another proof of this fact, also due
to Springer, which does not use non-abelian $H^2$, can be found in
[S2], III, \S 2.4. If the group $G$ is connected and linear, one can
therefore use Steinberg's theorem to conclude that each homogeneous
$G$-space has a $k$-point.

In this paper we prove a formally real analogue of Grothendieck's
theorem. That is, we assume that the ground field $k$ is perfect
and has {\it virtual\/} cohomological dimension one, $\vcd(k)\le1$.
This means that $k$ has some finite extension $K$ with $\cd(K)\le1$
(one can take $K=k(\sqrt{-1})$). 
% If $\vcd(k)<\cd(k)$ then $\char(k)=0$ and $k$ is perfect. 
Our main result specializes to Grothendieck's theorem if $\cd(k)\le1$, 
which is equivalent to 
% $\vcd(k)\le 1$ and  
$k$ being not formally real, i.e.\ having no orderings. However, even
in this case our proof is independent of the proof given in [Sp].
Typical examples of formally real fields with $\vcd(k)=1$ are
function fields of curves over $\Bbb R$, or the power series field
$\Bbb R (\!(T)\!)$; one can replace $\Bbb R$ by any real closed
field here ([S2], II, \S 3.3). Before we describe our results in 
more detail, let us briefly review previous work on formally real 
analogues of the theorems by Steinberg and Springer.

It was Colliot-Th\'el\`ene [CT] who proposed to study analogues of
the classical Hasse principle for the function field $k=\Bbb R(Y)$
of a smooth projective curve $Y$ over $\Bbb R$, in which the role
of the local places would be played by the completions $k_P$ of $k$
at the $\Bbb R$-points $P$ of $Y$. As observed by the second author
in [Sch], one can consider as local objects the real closures
$k_\xi$ of $k$ with respect to its orderings $\xi$. Apart from
having technical advantages, this point of view leads to stronger
results. The orderings of $k$ are the points of a compact, totally
disconnected topological space $\Omega_k=\operatorname{Sper}\,k$, the real
spectrum of $k$. The following Hasse principle was proved in [Sch]:
If $k$ is a perfect field with $\vcd(k)\le1$, $G$ is a connected
linear group over $k$ and $X$ is a $G$-torsor with a $k_\xi$-point
for each $\xi$ in a dense subset of $\Omega_k$, then $X$ has a
$k$-point. An equivalent way of expressing this is that the natural map
$$H^1(k,G)\longrightarrow\prod_\xi H^1(k_\xi,G)\eqno(*)$$
is injective, where in the product $\xi$ ranges over any dense
subset of the real spectrum $\Omega_k$. Note that this reduces to
Steinberg's theorem if $k$ has no orderings. An important point in
the proof is the notion of locally constant families of local
cohomology classes, technically realized through the construction
of a sheaf $\Cal H^1(G)$ on $\Omega_k$, and the description of the
precise image of $(*)$ in these terms.

The classical framework for Hasse principles is the context of
local and global fields, and in particular, of number fields. Here
one has exactly the same sort of Hasse principle for torsors if the
algebraic group $G$ is semisimple and simply connected. This is due
to work by Kneser and Harder in the 60s, and to more recent work of
Chernousov on $E_8$. Note that global or (nonarchimedean) local
fields $k$ have $\vcd(k)=2$. For torsors under semisimple, simply
connected classical groups (as well as $G_2$ and $F_4$), these
results have recently been generalized to all fields with $\vcd(k)=
2$, by Bayer-Fluckiger and Parimala [BP1], [BP2].

Assume again that ($k$ is perfect and) $\vcd(k)=1$. The paper [Sch] 
contains also the analogue of Springer's theorem (Thm.\ 6.5). Namely, 
if $G$ is an algebraic group over $k$ and $X$ is a homogeneous 
$G$-space such that over each real closure $k_\xi$ there exists a 
$G$-torsor which dominates $X$ (over $k_\xi$), then there exists a 
$G$-torsor over $k$ which dominates $X$ (over $k$). The proof used 
ideas of Springer's second proof in [S2]. Oddly, it needed the 
Feit-Thompson theorem on groups of odd order. The question of whether 
the assumption (of the existence of local dominating torsors) can be
relaxed to a dense subset of $\Omega_k$, instead of all of
$\Omega_k$, remained open and was raised in [Sch].

We shall now describe our main result. Assume that $k$ is perfect 
with $\vcd(k)\le1$, let $\ol G$ be an algebraic group over $k_s$ and 
$\kappa$ a $k$-kernel in $\ol G$. Our extension of Grothendieck's 
theorem takes the form of local-global principles for $H^2(k,\ol G, 
\kappa)$. Its main features are:
\item{a)} The set $H^2(k,\ol G,\kappa)$ contains a neutral element
if and only if $H^2(k_\xi,\ol G,\kappa)$ contains one for all $\xi$
in a dense subset of $\Omega_k$. In other words, $\kappa$ is
trivial if and only if $\kappa_\xi$ is trivial for all $\xi$ in a
dense subset of $\Omega_k$.
\item{b)} For $\alpha$, $\beta\in H^2(k,\ol G,\kappa)$, we have
that $\alpha=\beta$ if and only if $\alpha_\xi=\beta_\xi$ for all
$\xi$ in a dense subset of $\Omega_k$.
\item{c)} The element $\alpha$ is neutral if and only if $\alpha_
\xi$ is neutral for all $\xi$ in a dense subset of $\Omega_k$.
\par\n
In more technical terms, we construct a sheaf of sets $\Cal H^2(\ol
G,\kappa)$ on $\Omega_k$ which is locally constant and whose stalk
at $\xi$ is the finite set $H^2(k_\xi,\ol G,\kappa)$. Further we
construct a subsheaf $\Cal N^2(\ol G,\kappa)$ of $\Cal H^2(\ol G,
\kappa)$ which is again locally constant and whose stalk at $\xi$
is the subset $N^2(k_\xi,\ol G,\kappa)$ of $H^2(k_\xi,\ol G,
\kappa)$. Finally we show that the natural map $H^2(k,\ol G,\kappa)
\to\Gamma(\Omega_k,\Cal H^2(\ol G,\kappa))$ is bijective, and its
restriction to $N^2(k,\ol G,\kappa)$ bijects to $\Gamma(\Omega_k,
\Cal N^2(\ol G,\kappa))$.

The most difficult step is the proof of c). We also present an
alternative approach, which is technically easier but applies only
when
% $\ch(k)=0$ and
$\ol G$ is connected and linear. It is based on Borovoi's elegant
technique of hypercohomological abelianization of the 
non-commutative $H^2$ ([B1], [B2]). Both approaches use the technique 
of sheafified $H^1$ from [Sch].

As an application of our extension of Grothendieck's theorem, we
give a new and simpler proof to the formally real analogue of
Springer's theorem ([Sch], Thm.\ 6.5). It is the fact that we
know the sheaves $\Cal H^2$ and $\Cal N^2$ to be locally constant,
which actually allows us to work with just a dense subset of
$\Omega_k$. In particular, we obtain an affirmative answer to the
question mentioned above, of whether domination of a homogeneous
space locally over a dense subset of $\Omega_k$ suffices to
conclude the existence of a dominating $k$-torsor.

The paper is organized as follows. In Section~1 we introduce
kernels and non-abelian $H^2$ for algebraic groups, as well as the
subset $N^2$ of neutral elements. Section~2 constructs the locally
constant sheaves $\Cal H^2$ and $\Cal N^2$ and studies their
properties. The formally real analogues of Grothendieck's theorem
are proven in Section~3. Section~4 presents the alternative
approach based on Borovoi's abelianization technique. Finally,
applications to local-global principles for homogeneous spaces are
discussed in Section 5.

\smallskip
We would like to thank J.-P.~Serre for his critical remarks which 
helped to improve the exposition. The first author thanks also
P.~Deligne for illuminating correspondence, TIFR, Bombay, and Tel-Aviv 
University -- in particular J.~Bernstein, M.~Borovoi, M.~Jarden -- for 
hospitality and interest, and the NSF for the grant INT-9603014. The
first and second authors thank NATO for the grant CRG-970133. The third 
author gratefully acknowledges the support of Alexander von Humboldt 
foundation and the hospitality of Universit\"at Regensburg.
\bigbreak

\n
{\bf Notations and conventions}
\par\n
Let $k$ be a field with separable closure $k_s$ and absolute Galois
group $\Gamma_k=\Gal(k_s/k)$. The cohomological dimension $\cd(k)$ 
of $k$ is the largest integer $n$ for which there is a finite discrete
$\Gamma_k$-module $A$ with $H^n(\Gamma_k,A)\ne0$ (resp.\ is
$\infty$ if no largest such $n$ exists). The virtual cohomological
dimension $\vcd(k)$ is the common cohomological dimension of all
sufficiently large finite separable extensions of $k$; 
equivalently, $\vcd(k)=\cd k(\sqrt{-1})$. Note that $\vcd
(k)<\cd(k)$ can happen only if $k$ has an ordering, in particular,
only if $\ch(k)=0$ ([S2], II, Prop. 4.1).
Generally we do not assume the base field $k$ to be perfect, but 
our main results will need this hypothesis. 

The real spectrum of $k$, which we denote here by $\Omega_k$ and
which is often written $\Sper k$, is the topological space of all
orderings $\xi$ of $k$. Its topology is generated by the subsets
$\{\xi\in\Omega_k\colon$ $a>0$ at~$\xi\}$ for $a\in k$. This is a
boolean (=~compact and totally disconnected) topological space.
Given $\xi\in\Omega_k$, one denotes the real closure of $k$ at
$\xi$ by $k_\xi$. Note that $\Omega_k$ is naturally homeomorphic to
the topological quotient space of $\Inv(\Gamma_k)$ (the space of
elements of $\Gamma_k$ of order two) modulo conjugation by
$\Gamma_k$. See e.g.\ [Scha], ch.~3, \S5, for some basic information
on $\Omega_k$.

Throughout the paper, by an algebraic group over $k$ we mean a {\it
smooth\/} group scheme of finite type over $\Spec k$. If $\ch(k)=
0$, the smoothness assumption holds automatically ([DG], II, \S6,
no.~1). An algebraic group over $k$ which is absolutely reduced 
(reduced over an algebraic closure $\ol k$ of $k$) is smooth over $k$ 
([DG], II, \S 5, 2.1(v)).
\bigskip
\subhead 1. Non-commutative $H^2$ for algebraic groups\endsubhead

\n In this section we define kernels and noncommutative $H^2$ for (not
necessarily linear) algebraic groups. We will follow mainly Borovoi
[B1], who rewrote part of Springer [Sp], but only for linear groups
and fields of characteristic zero. For a general account in terms
of gerbes see Giraud [G], Deligne-Milne [DM], Breen [Br].

\lab(1.1) Let $k$ be a field and $k_s$ a fixed separable closure of
$k$. Denote the profinite group $\Gal(k_s/k)$ by $\Gamma_k$, or
just by $\Gamma$. Given $s\in\Gamma_k$, let $s^*$ denote the
morphism $\Spec k_s\to\Spec k_s$ induced by $s$. Note that $(st)^*=
t^*s^*$.

\lab(1.2) Let $\ol G$ be an algebraic group over $k_s$. Denote by
$\Aut(\ol G)$ the group of automorphisms of $\ol G$ (where $\ol G$ 
is considered as a group scheme over $k_s$). Given $s\in\Gamma_k$, let 
$s_*\ol G$ denote the base change of $p\colon\ol G\to\Spec k_s$ by $s^*$. 
%In other words, if $p\colon\ol G\to\Spec k_s$ is the structural morphism 
%of $\ol G$ then $s_*\ol G$ is the group scheme $\ol G$ with structural 
%morphism $(s^*)^{-1}\comp p$ to $\Spec k_s$. 
Then $s_*\ol G$ is another group scheme over $\Spec k_s$, which as a
$k_s$-scheme is isomorphic to $(s^*)^{-1}\comp p\colon$ $\ol G\to\Spec k_s$.

An {\it $s$-semilinear automorphism\/}
of $\ol G$ is by definition an isomorphism of algebraic groups over
$k_s$ from $s_*\ol G$ to $\ol G$. A {\it $k$-semilinear automorphism\/}
$\varphi$ of $\ol G$ will mean an $s$-semilinear automorphism of $\ol G$
for some $s\in\Gamma_k$. % It is easy to see 
Note that $s$ is uniquely determined by $\varphi$, since $s^\ast=
p\comp\varphi\comp e$, where $e\colon\Spec k_s\to\ol G$ is the identity point. 

The set of $k$-semilinear automorphisms of $\ol G$ forms a group which we
denote by $\SAut(\ol G/k)$. If $\varphi\colon s_*\ol G\to\ol G$ (resp.\ 
$\psi\colon t_*\ol G\to\ol G$) is an $s$-semilinear (resp.\ $t$-semilinear)
automorphism of $\ol G$, then the product $\psi\cdot\varphi$ is by 
definition the $t\comp s$-semilinear automorphism $\psi\comp t_*\varphi\colon$ 
$t_*s_*\ol G\to\ol G$ of $\ol G$. 
Sending an $s$-semilinear automorphism to $s$ defines an exact sequence
$$1\to\Aut(\ol G)\to\SAut(\ol G/k)\to\Gamma_k.\eqno(1)$$
The last map need not be surjective in general.

\lab(1.3) There is a natural action of the group $\SAut(\ol G/k)$
on the group $\ol G(k_s)$, given as follows. For $s\in\Gamma_k$ let
$\beta_s\colon\>\ol G(k_s)\to(s_*\ol G)(k_s)$ be defined by $\beta_
s(x):=x\comp s^*$. Then $\beta_s$ is an isomorphism of groups. If
$\varphi\colon s_*\ol G\to\ol G$ is an $s$-semilinear automorphism, define
$\varphi_*\colon\>\ol G(k_s)\to\ol G(k_s)$ by composing $\varphi$
with $\beta_s$, i.e.\ $\varphi_*(x):=\varphi\comp x\comp s^*$ for
$x\in\ol G(k_s)$. Then $\varphi_*$ is an automorphism of the group
$\ol G(k_s)$; moreover $(\psi\varphi)_*=\psi_*\comp\varphi_*$ if
$\psi$ is another semilinear automorphism of $\ol G$.

Thus we have defined a group homomorphism $\SAut(\ol G/k)\to\Aut\ol
G(k_s)$, $\varphi\mapsto\varphi_*$. In general this homomorphism
need not be injective, e.g.\ if $\ol G$ is a finite constant group
scheme. Usually we will simply write $\varphi(x)$ instead of
$\varphi_*(x)$.

Observe that $\SAut(\ol G/k)$ acts also on $k_s[\ol G]:=\Gamma(\ol
G,\Cal O_{\ol G})$. In this case the action is given by $\varphi^*
(a):=a\comp\varphi$ ($\varphi\in\SAut(\ol G/k)$, $a\in k_s [\ol G]$; 
here we regard $\varphi$ as a morphism of schemes $\ol G\to\ol G$ 
which satisfies $p\comp\varphi=(s^*)^{-1}\comp p$).
% the semilinear automorphisms of $\ol G$ are scheme morphisms $\ol G\to
% \ol G$.

\lab(1.4) A {\it$k$-form\/} of $\ol G$ is an algebraic group $G$
over $k$ together with an isomorphism $\ol G\cong G\times_kk_s$ of
algebraic groups over $k_s$. If a $k$-form of $\ol G$ is fixed, the
group $\ol G$ is also said to be defined over $k$. Any $k$-form of
$\ol G$ defines a splitting $\Gamma_k\to\SAut(\ol G/k)$ of $(1)$,
by $s\mapsto\operatorname{id}\times_k(s^{-1})^*$. For a converse see (1.15)
below.

\lab(1.5) Given $x\in\ol G(k_s)$ we write $\int(x)$ for the inner
automorphism $y\mapsto xyx^{-1}$ of $\ol G$. Let $\Int(\ol G)$ be
the subgroup of $\Aut(\ol G)$ consisting of the $\int(x)$, $x\in\ol
G(k_s)$. Thus $\Int(\ol G)=\ol G(k_s)/\ol Z(k_s)$ where $\ol Z$ is
the center of $\ol G$. The subgroup $\Int(\ol G)$ is normal in
$\SAut(\ol G/k)$. Let
$$\Out(\ol G)\>:=\>\Aut(\ol G)/\Int(\ol G)$$
and
$$\SOut(\ol G/k)\>:=\>\SAut(\ol G/k)\big/\Int(\ol G).$$
Taking $(1)$ modulo $\Int(\ol G)$ we get the exact sequence
$$1\to\Out(\ol G)\to\SOut(\ol G/k)\to\Gamma_k.\eqno(2)$$ \bigbreak

\n{\bf(1.6) Definition.}\enspace We equip $\SAut(\ol G/k)$
with the weak topology with respect to the family of evaluation
maps $\operatorname{ev}_x\colon\>\SAut(\ol G/k)\to\ol G(k_s)$, $\varphi
\mapsto\varphi(x)$ ($x\in\ol G(k_s)$, see (1.3)),
% ie. the topology induced by the natural inclusion of $\SAut$ into
% the direct product $\prod_{\ol G(k_s)}\ol G(k_s)$,
where $\ol G(k_s)$ is given the discrete topology. A map $t\mapsto
f_t$ from a topological space $T$ to $\SAut(\ol G/k)$ will be
called {\it weakly continuous\/} if it is continuous with respect
to this topology. This is the case if and only if for every $x\in
\ol G(k_s)$ the map $T\to\ol G(k_s)$, $t\mapsto f_t(x)$ is
continuous (=~locally constant), if and only if the natural map $T
\times\ol G(k_s)\to\ol G(k_s)$ is continuous. \bigbreak

\n{\bf(1.7) Proposition.\enspace\it Let $\ol G$ be an
algebraic group over $k_s$ and $f\colon\>\Gamma_k\to\SAut(\ol G/k)$
a set-theoretic section of\/ $(1)$. Let $K\supset k$ be a finite
Galois extension such that there is a $K$-form $\widetilde G$ of
$\ol G$. Let $\sigma \colon\>\Gamma_K\to\SAut(\ol G/K)\subset\SAut
(\ol G/k)$ be the splitting of\/ $(1)$ associated with $\widetilde
G$. Consider the following conditions:
\item{\rm(i)} $f$ is weakly continuous (cf.\ (1.6));
\item{\rm(ii)} for every $s\in\Gamma_k$ the map
$$\Gamma_K\to\Aut(\ol G),\quad t\mapsto\sigma_t^{-1}f_s^{-1}f_{st}
\eqno(3)$$
is locally constant.
\par\n
Then\/ {\rm(ii)} implies\/ {\rm(i)}, and the 
converse holds if\/$\ch(k)=0$. 

Moreover, if $\ol G$ is linear and\/ $\ch(k)$ is
arbitrary, {\rm(ii)} is equivalent to
\item{\rm(iii)} for every $a\in k_s[\ol G]$ the map
$$\Gamma_k\to k_s[\ol G],\quad s\mapsto f_s^*(a)\eqno(4)$$
is locally constant. \bigbreak}

\Proof Let $K$, $\widetilde G$ and $\sigma$ be as in the
proposition. Fix $s\in\Gamma_k$ and write $\varphi_t:=\sigma_t^{-1}
f_s^{-1}f_{st}$ for $t\in\Gamma_K$, so that $f_{st}=f_s\sigma_t
\varphi_t$. If the map $\Gamma_K\to\Aut(\ol G)$, $t
\mapsto\varphi_t$ is locally constant then obviously the map
$\Gamma_K\to\ol G(k_s)$, $t\mapsto f_{st}(x)$ is locally constant
for every $x\in\ol G(k_s)$. From this one concludes that (ii)
implies (i).

Assume that $\ol G$ is linear. Since $k_s[\ol G]$ is a finitely
generated $k_s$-algebra, the map $(3)$ is locally constant if and
only if for every $a\in k_s[\ol G]$ the map $\Gamma_K\to k_s[\ol
G]$, $t\mapsto\varphi_t^*(a)$ is locally constant. On the other
hand, (iii) can be reformulated as saying that for every $s\in
\Gamma_k$ and every $a\in k_s[\ol G]$ the map $\Gamma_K\to k_s[\ol
G]$, $t\mapsto f_{st}^*(a)$ is locally constant. Since $f_{st}=f_s
\sigma_t\varphi_t$, it is clear that (ii) and (iii) are equivalent.

Finally assume $\ch(k)=0$. Then the implication (i) $\Rightarrow$
(ii) follows from the next lemma, from which it follows that in
order to verify that $t\mapsto\varphi(t)$ is locally constant, we
only need to check that $t\mapsto\varphi_t(x)$ is, for finitely
many $x\in\ol G(k_s)$: \bigbreak

\n{\bf(1.8) Lemma.\enspace\it Let $\ol k$ be an
algebraically closed field of characteristic zero. Then for every
algebraic group $\ol G$ over $\ol k$ there exists a finitely
generated subgroup $S$ of $\ol G(\ol k)$ which is Zariski dense in
$\ol G$. \medbreak}

\Proof We can assume that $\ol G$ is connected. Let $\ol H$ be the
(unique) maximal element among the identity connected components of
Zariski closures of finitely generated subgroups of $\ol G(k_s)$.
If $\ol H\ne\ol G$ then $\ol G/\ol H$ is a non-trivial connected
algebraic group for which the group $(\ol G/\ol H)(\ol k)$ is
locally finite (i.e.\ every finitely generated subgroup is finite).
But every non-trivial connected algebraic group over $\ol k$
contains an element of infinite order. This is obvious for $\Bbb
G_a$ and $\Bbb G_m$, and is known for abelian varieties (e.g.\
[FJ], Theorem 10.1), from which the general case follows. Therefore
we must have $\ol H=\ol G$. \endproof

\n{\bf(1.9) Remarks.}
\par\n
1.\enspace Lemma (1.8) is also valid if $\ol k$ is a separably
closed field of positive characteristic which is not the algebraic
closure of a finite field, and if $\ol G$ has no nontrivial
unipotent quotient groups (same proof). Therefore the implication
(i) $\Rightarrow$ (ii) of Proposition (1.7) holds in this case as
well.
\par\n
2.\enspace In characteristic zero, Lemma (1.8) implies that the
weak topology on $\Aut(\ol G)$ (cf.\ (1.6)) is discrete. In
positive characteristics this is in general not true. For example,
the group $\ol G=\Bbb G_a\times\Bbb G_a$ has automorphisms of the
form $(x,y)\mapsto\bigl(x+P(y),\,y\bigr)$ where $P$ is any additive
polynomial. These automorphisms form a subgroup of $\Aut(\ol G)$
which is not discrete.

From this example it follows easily that in (1.7), (i) does not
always imply (ii) and (iii) if $\ch(k)=p>0$. \bigbreak

\n{\bf(1.10) Definition.}\enspace A section $f$ of $(1)$
will be called {\it continuous\/} if it satisfies condition (ii) of
Proposition (1.7). If $f$ is continuous then it is weakly
continuous, the converse being true if $\ch(k)=0$ (1.7). \bigbreak

\n{\bf(1.11) Definition.}\enspace A {\it$k$-kernel in $\ol
G$} is a group homomorphism $\kappa\colon\>\Gamma_k\to\SOut(\ol
G/k)$ which splits $(2)$ and lifts to a continuous {\it map\/}
$f\colon\>\Gamma_k\to\SAut(\ol G/k)$. A pair $L=(\ol G, \kappa)$
consisting of an algebraic group $\ol G$ over $k_s$ and a
$k$-kernel $\kappa$ in $\ol G$ will simply be called a {\it
$k$-kernel}. Other terms would be $k$-lien or $k$-band [G], [DM],
[Br]. \bigskip

\n{\bf(1.12) Lemma.\enspace\it The weak topology on\/ $\Int
(\ol G)$ is discrete. \medbreak}

\Proof We have to show that there is a finitely generated subgroup
$S$ of $\ol G(k_s)$ whose centralizer is the center of $\ol G
(k_s)$. Since the Zariski topology is noetherian, there is a
minimal group among all centralizers of finitely generated
subgroups of $\ol G(k_s)$. This group is necessarily the center of
$\ol G(k_s)$. \endproof

\n{\bf(1.13) Corollary.\enspace\it Let $\kappa$ be a
splitting of\/ $(2)$, and let $f$, $f'$ be two weakly continuous
maps from\/ $\Gamma_k$ to $\SAut(\ol G/k)$ each of which lift
$\kappa$. If $f$ is continuous then so is $f'$. \medbreak}

\Proof This follows from Lemma (1.12). \endproof

\n{\bf(1.14) Remarks.}\enspace We add a few comments here on
other definitions of kernels of algebraic groups found in the
literature. Let $\ol G$ be an algebraic group over $k_s$.

Springer ([Sp], p.~176, bottom) defines more generally
$K/k$-kernels for any Galois extension $K/k$. We consider his
definition only in the case $K=k_s$. In [Sp], 
%\footnote{$^2$}
%  {modulo an amazing number of misprints within only a few lines}
a $k_s/k$-kernel in $\ol G$ is a group homomorphism $\Gamma_k\to
\Out\ol G(k_s)$, $s\mapsto\kappa(s)$, such that $\kappa(s)$ is
induced by an $s$-semilinear automorphism of $\ol G$ ($s\in\Gamma_k$) and
such that $\kappa$ is a locally constant map. Assume --- as is done
in [Sp], end of p.~176 --- that the canonical map $\SAut(\ol G/k)
\to\Aut\ol G(k_s)$ is injective, cf.\ (1.3).
% (In general it is not, but in most cases it is...)
Then $\SOut(\ol G/k)\to\Out\ol G(k_s)$ is also injective, so
$\kappa$ can be viewed as a locally constant homomorphism
$\Gamma_k\to\SOut(\ol G/ k)$. But then {\it a fortiori\/} its
composition with the canonical projection $\SOut(\ol
G/k)\to\Gamma_k$ is locally constant. By definition, this
composition is the identity of $\Gamma_k$, which implies that
$\Gamma_k$ is discrete.

Borovoi takes the problem of continuity into account. He only
considers the case where $\ch(k)=0$ and the group $\ol G$ is linear
([B1], 1.3(b)). His definition of a $k$-kernel $\kappa$ in $\ol G$
is the same as our Definition (1.11) (in fact, it inspired our
definition), except that his notion of continuity, for maps $f
\colon\Gamma_k\to\SAut(\ol G/k)$ which are sections of $(1)$, is
different. His condition requires that for any $a\in k_s[\ol G]$
the stabilizer of $a$ in $\Gamma_k$ be open. Since $f$ is not
assumed to be a homomorphism, this stabilizer is only a {\it
subset\/}, not a subgroup, in general. Examples show that his
condition does not imply condition (iii) in Proposition (1.7),
which (under his hypotheses) is equivalent to our notion of
continuity (and, under his hypotheses, also to weak continuity).
Therefore Borovoi's definition seems too weak in general. \bigbreak

\n{\bf(1.15) Remark.}\enspace The splitting of $(1)$ defined
by a $k$-form $G$ of $\ol G$ is continuous. Hence reading this
splitting modulo $\Int(\ol G)$ gives a $k$-kernel in $\ol G$. We
will denote it by $\kappa_G$. Conversely, a continuous splitting of
$(1)$ defines a unique $k$-form of $\ol G$. (See [BS], Lemme 2.12;
condition (b) holds by our definition of continuity, while (c)
holds by [S1], ch.~V no.~20, Cor.~2.) If $\ch(k)=0$, therefore,
every weakly continuous splitting of $(1)$ defines a unique
$k$-form of $\ol G$. \bigbreak

% \n{\bf(1.XX) Corollary.\enspace\it Any two continuous
% sections $f$, $f'$ of $(1)$ differ by a map $\Gamma\to\Aut(\ol
% G)$ which is locally constant. \endproof}(in char=0)

\lab(1.16) If $L=(\ol G,\kappa)$ is a $k$-kernel and $\ol Z$ is the
center of $\ol G$ then $\kappa$ induces a $k$-kernel in $\ol Z$,
that is, $\kappa$ defines a $k$-form $Z$ of $\ol Z$. We call $Z$
(which is a commutative algebraic group over $k$) the {\it
center\/} of the $k$-kernel $L$. \medbreak

Let now $L=(\ol G,\kappa)$ be a $k$-kernel. Recall the definition
of the cohomology set $H^2(k,L)$: \bigbreak

\n{\bf(1.17) Definition.\enspace} A {\it$2$-cocycle with
coefficients in $L$} is a pair $(f,g)$ of maps
$$f\colon\>\Gamma\to\SAut(\ol G/k),\ s\mapsto f_s,\quad{\text and}
\quad g\colon\>\Gamma\times\Gamma\to\ol G(k_s),\ (s,t)\mapsto g_{s,
t},$$
such that
\item{(1)} $f$ is continuous (1.10), and $f$~mod~$\Int(\ol G)=
\kappa$;
\item{(2)} $g\colon\>(s,t)\mapsto g_{s,t}$ is continuous (= locally
constant), and for $s,t,u\in\Gamma$ one has
$$f_s\comp f_t\>=\>\int(g_{s,t})\comp f_{st}\quad{\text and}\quad
f_s(g_{t,u})\cdot g_{s,tu}\>=\>g_{s,t}\cdot g_{st,u}.\eqno(5)$$
\par\n
Let $Z^2(k,L)$ denote the set of these 2-cocycles. Two cocycles
$(f,g)$ and $(f',g')$ are called {\it equivalent\/} if there is a
continuous (= locally constant) map $h\colon\>\Gamma\to\ol G(k_s)$
such that
$$f'_s\>=\>\int(h_s)\comp f_s\quad{\text and}\quad g'_{s,t}\>=\>h_s
\cdot f_s(h_t)\cdot g_{s,t}\cdot h_{st}^{-1}\eqno(6)$$
for all $s,t\in\Gamma$. The cohomology set $H^2(k,L)$ is defined to
be the set of equivalence classes in $Z^2(k,L)$.
If $\ol G$ is commutative then $H^2(k,L)$ is the usual second
(Galois) cohomology group $H^2(k,G)$, where $G$ is the $k$-form of
$\ol G$ defined by $\kappa$.

\lab(1.18) The description of cocycles above follows [Sp], rather
than [B1] who takes $(f,g^{-1})$ instead of $(f,g)$. An alternative
and useful (see (3.4) and (3.5) below) description of $H^2(k,L)$ in
terms of group extensions is recalled next (compare [Sp]). Consider
extensions
$$1\to\ol G(k_s)\rightmap iE\rightmap\pi\Gamma\to1\eqno(7)$$
of topological groups, where $\ol G(k_s)$ (resp.\ $\Gamma=\Gamma_
k$) carries the discrete (resp.\ the natural profinite) topology,
and $i$ and $\pi$ are open onto their respective images. Two such
extensions $E$ and $E'$ are called {\it equivalent\/} if there is
an isomorphism $E\to E'$ of topological groups which induces the
identity on $\ol G(k_s)$ and on $\Gamma$. \bigbreak

\n{\bf(1.19) Lemma.\enspace\it The set $H^2(k,L)$ is in
natural bijection with the set of equivalence classes of
extensions\/ $(7)$ for which the induced homomorphism from\/
$\Gamma$ to $\Out\ol G(k_s)$ coincides with the composite
homomorphism
$$\Gamma\rightmap\kappa\SOut(\ol G/k)\to\Out\ol G(k_s).\eqno(8)$$
The second map in $(8)$ is induced by the canonical map $\SAut(\ol
G/k)\to\Aut\ol G(k_s)$ (1.3). \medbreak}

\Proof Given such an extension $(7)$, choose a continuous section
$\Gamma\to E$, $s\mapsto z_s$. Such sections exist since $\ol G(k_
s)$ is a discrete subgroup of $E$. For any $s\in\Gamma$ there is a
unique $s$-semilinear automorphism $f_s$ of $\ol G$ with $f_s(x)=z_sxz_s
^{-1}$ for every $x\in\ol G(k_s)$. Put $g_{s,t}:=z_sz_tz_{st}^{-1}$
for $s,t\in\Gamma$. The pair $(f,g)$ is a 2-cocycle with values in
$L$. Indeed, it is obvious that $f$ is weakly continuous. Since
$\kappa$ has a continuous lift by hypothesis, $f$ is in fact
continuous by Corollary (1.13). Choosing a different continuous
section $s\mapsto z'_s$ yields an equivalent cocycle.

Conversely, to a cocycle $(f,g)\in Z^2(k,L)$ one associates an
extension $(7)$ by putting
$$\hbox{$E\>:=\>\ol G(k_s)\times\Gamma$ \ with multiplication rule
\ $(x,s)\cdot(y,t):=\bigl(x\,f_s(y)\,g_{s,t},\,st\bigr)$.}$$
The group $E$ is given the product topology, and $i$, $\pi$ are
defined in the obvious way. The equivalence class of this extension
depends only on the cohomology class of $(f,g)$, and the two
processes are inverses of each other. \endproof

% THE FOLLOWING IS TRUE ONLY IN char = 0! HOWEVER WE CAN GET BY
% WITHOUT THIS ARGUMENT LATER.
%
% \lab(1.YY) The above argument shows also the following: A
% topological extension $(7)$ belongs to some $k$-kernel if (and
% only if) for every $s\in\Gamma$ the outer automorphism of $\ol
% G(k_s)$ associated to $s$ by the extension is induced by some
% $s$-semilinear automorphism of $\ol G$. Indeed, assuming the latter,
% choose a continuous section $s\mapsto z_s$ of $(7)$. Then, as
% above, there is a unique $s$-semilinear automorphism $f_s$ which acts as
% $\int(z_s)$ on $\ol G(k_s)$. The map $s\mapsto f_s$ is
% continuous, and $\kappa(s):=\ol{f_s}\in\SOut(\ol G/k)$ is a
% kernel (in fact, the unique kernel) to which $(7)$ belongs.

\lab(1.20) Given a $k$-kernel $L=(\ol G,\kappa)$, the set $H^2(k,
L)$ may be empty. One has the following well-known criterion ([M],
IV, \S8; [G], VI, \S2). Let $Z$ be the center of $\ol G$, considered
as an algebraic group over $k$ via $\kappa$. With $L$ one
associates a class $\obs(L)$ in the cohomology group $H^3(k,Z)$ as
follows: 
% Mac Lane writes $\Obs(L)$, Giraud writes $c(L)$

Let $s\mapsto f_s$ be a continuous map $\Gamma\to\SAut(\ol G/k)$
which lifts $\kappa$, cf.\ (1.11). The map $(s,t)\mapsto f_s
f_tf_{st}^{-1}$ from $\Gamma\times\Gamma$ to $\Int(\ol G)$ is
locally constant by Lemma (1.12), hence there exists a locally
constant map $\Gamma\times\Gamma\to\ol G(k_s)$, $(s,t)\mapsto g_
{s,t}$ such that
$$f_s\comp f_t\>=\>\int(g_{s,t})\comp f_{st}.$$
Let $z\colon\>\Gamma\times\Gamma\times\Gamma\to\ol Z(k_s)$ be the
locally constant map determined by
$$f_s(g_{t,u})\cdot g_{s,tu}\>=\>z_{s,t,u}\cdot g_{s,t}\cdot g_{st,
u}.$$
Then $z\in Z^3(k,Z)$, and the class of $z$ in $H^3(k,Z)$ is
independent of the choices made. Denote it by $\obs(L)$. One has:
\bigbreak

\n{\bf(1.21) Proposition.\enspace\it The set $H^2(k,L)$ is
non-empty if and only if\/ $\obs(L)=0$ in $H^3(k,Z)$. \medbreak}

\Proof See [M], IV, Thm.\ 8.7. \endproof

This yields the following local-global principle for non-emptiness
of $H^2$: \bigbreak

\n{\bf(1.22) Corollary.\enspace\it Assume\/ $\vcd(k)\le2$,
and let $L=(\ol G,\kappa)$ be a $k$-kernel. Then $H^2(k,L)\ne
\emptyset$ if and only if $H^2(k_\xi,L)\ne\emptyset$ for all $\xi$
in a dense subset of\/ $\Omega_k$. 
\medbreak}

For the proof we need a straightforward (partial) generalization of
[Sch], Thm.\ 3.1, which we formulate in a greater generality than
actually required here. The proof is analogous to {\it loc.cit.}
Recall that sheaves $\Cal H^i(A)$ on $\Omega_k$ are attached in [Sch], 
(2.8), to any commutative algebraic group $A$ over $k$ and $i\ge 1$. 
They are shown to be locally constant in [Sch], (2.13(a)), and to have 
stalks $H^i(k_\xi,A)$ in [Sch], (2.9(b)).
\bigbreak

\n{\bf(1.23) Lemma.\enspace\it If $k$ is perfect with\/
$\vcd(k) \le d$ and $A$ is a commutative algebraic group over $k$,
then the map $H^n(k,A)\to\Gamma\bigl(\Omega_k,\,\Cal H^n(A)\bigr)$
is bijective for $n>d$ and surjective for $n=d$. \endproof}

\n
{\it Proof\/} of Corollary (1.22). Assume $H^2(k_\xi,L)\ne
\emptyset$ for all $\xi$ in a dense subset of $\Omega_k$. Let $Z$
be the center of $L$. Writing $L_\xi$ for the restriction of $L$ to
$k_\xi$, $\obs(L)$ maps to $\obs(L_\xi)$ under the restriction map
$H^3(k,Z)\to H^3(k_\xi,Z)$. By (1.21), and by the assumption, $\obs
(L)$ lies in the kernel of $H^3(k,Z)\to\Gamma\bigl(\Omega_k,\,\Cal
H^3(Z)\bigr)$. By Lemma (1.23) this kernel is zero. \endproof

\lab(1.24) Let $L$ be a $k$-kernel with center $Z$. There is a
natural action of the abelian group $H^2(k,Z)$ on the set $H^2(k,
L)$, which is free and transitive provided $H^2(k,L)\ne\emptyset$.
In fact, the abelian group $Z^2(k,Z)$ acts on the set $Z^2(k,L)$ by
$z\cdot(f,g)=(f,zg)$, and this action descends to an action of $H^2
(k,Z)$ on $H^2(k,L)$.

For completeness we also mention the interpretation of the action
in terms of group extensions: Given an extension $(7)$ whose class
in $H^2(k,L)$ is $\alpha$, and an extension
$$1\to\ol Z(k_s)\to B\to\Gamma\to1$$
whose class in $H^2(k,Z)$ is $\zeta$, form the extension
$$1\to\ol G(k_s)\to(B\times_\Gamma E)/D\to\Gamma\to1\eqno(9)$$
where $D$ is the subgroup $\{(z,z^{-1})\colon$ $z\in\ol Z(k_s)\}$
of the fiber product. Then $(9)$ has a class in $H^2(k,L)$, and
this class is $\zeta\cdot\alpha$.

Observe that the action is compatible with extension of the base
field, i.e.\ with restriction in cohomology. \bigbreak

\n{\bf(1.25) Definition.}\enspace A 2-cocycle $(f,g)\in Z^2
(k,L)$ is called {\it neutral\/} if $g_{s,t}=1$ for all $s,t$. An
element $\alpha\in H^2(k,L)$ is called {\it neutral\/} if it can be
represented by a neutral cocycle. The subset of $H^2(k,L)$
consisting of the neutral elements is denoted $N^2(k,L)$. The
kernel $L=(\ol G,\kappa)$ is called {\it trivial\/} if $N^2(k,L)$
is nonempty. \bigbreak

\lab(1.26) In terms of group extensions the neutral elements have
the following description: Given an extension $(7)$ compatible with the
kernel $\kappa$, its class in $H^2(k,L)$ is neutral if and only if
the extension splits by a continuous homomorphic section $\Gamma\to E$.

\lab(1.27) The kernel $L=(\ol G,\kappa)$ is trivial (i.e.\ $N^2(k,
L)\ne\emptyset$, (1.25)) if and only if the homomorphism $\kappa
\colon\>\Gamma\to\SOut(\ol G/k)$ lifts to a continuous homomorphism
$\Gamma\to\SAut(\ol G/k)$.
% Suppose there exists a continuous homomorphic lift $f$ of
% $\kappa$. Then we define $g_{s,t}=1$ for all $s,t\in\Gamma$ and
% $(f,1)$ is a neutral element. Conversely, suppose $(f,1)$ is a
% neutral element. Then $f_sf_t=f_{st}$ and $s\mapsto f_s$ is a
% continuous homomorphism.
By (1.15) it is equivalent to saying that $\kappa$ belongs to some
$k$-form $G$ of $\ol G$. Observe that $N^2(k,L)$ can be empty even
if $H^2(k,L)$ is non-empty (see also (3.6) below).
\bigskip
\subhead 2. The sheaves $\Cal H^2$ and $\Cal N^2$\endsubhead
%\beginsect2.\ The sheaves $\Cal H^2$ and $\Cal N^2$

\n Given an algebraic group $G$ over $k$, a natural sheaf $\Cal H^1
(G)$ of pointed sets on the space $\Omega_k$ of orderings of $k$,
was defined in [Sch]. It is locally constant and its stalk at $\xi$
is $H^1(k_\xi,G)$. %Also, 
As noted in (1.22), if $G$ is commutative, locally constant
sheaves $\Cal H^n(G)$ of abelian groups were defined for all
$n\ge1$; they have the analogous property of the stalks.

We now want to sheafify the non-commutative cohomology sets $H^2$
and $N^2$ in a similar way. We frequently present sheaves on
boolean (i.e.\ compact and totally disconnected) spaces (such as
$\Omega_k$) just by giving the sections over clopen subsets. This
is justified by [Sch], Appendix C.1.

\lab(2.1) Let $L=(\bar G,\kappa)$ be a $k$-kernel, with center $Z$.
The cohomology set $H^2(E,L)$ is defined for every finite separable
extension $E$ of $k$ and, more generally, for every finite \'etale
$k$-algebra $E$. In the latter case, if $E=K_1\times\cdots\times
K_r$ with $K_i/k$ finite separable field extensions, then $H^2(E,L)
=H^2(K_1,L)\times\cdots\times H^2(K_r,L)$. Therefore we can imitate
the definition in 2.3 of [Sch]: \bigbreak

\n{\bf(2.2) Definition.}\enspace The sheaf (of sets) $\Cal
H^2(L)$ on $\Omega_k$ is defined by
$$U\quad\longmapsto\quad\dirlim{(E,s)\in J_U}H^2(E,L),$$
for $U\subset\Omega_k$ clopen. Here $J_U$ is the category of pairs
$(E,s)$ where $E$ is an \'etale $k$-algebra and $s\colon U\to\Omega_E$ 
is a section of the restriction map $\Omega_E\to\Omega_k$ over
$U$ (see [Sch], Sect. 2.3).

Observe that one has a canonical map 
$H^2(k,L)\to\Gamma\bigl(\Omega_k,\,\Cal H^2(L)\bigr)$. \medbreak

Next we identify the stalks of the sheaves $\Cal H^2(L)$. For this
we need \bigbreak

\n{\bf(2.3) Lemma.\enspace\it If $K\supset k$ is a separable
algebraic extension then $H^2(K,L)$ is the direct limit of the sets
$H^2(F,L)$, where $F$ ranges over the finite subextensions $F
\supset k$ of $K\supset k$. \medbreak}

\Proof From Proposition (1.21) it follows that the lemma holds if
$H^2(F,L)=\emptyset$ for all these $F$. Hence we may assume $H^2(F,
L)\ne\emptyset$ for some $F$. Fixing an element $\alpha\in H^2(F,L)$ 
we get a bijection $H^2(F',Z)\isomap H^2(F',L)$ for every $F'
\supset F$, and these bijections are compatible with restriction
(1.24). So the assertion follows from the corresponding fact for
abelian cohomology of $Z$, which is well known. \endproof

\n{\bf(2.4) Corollary.\enspace\it The sheaf $\Cal H^2(L)$ is
locally constant with finite stalks. The stalk at $\xi\in\Omega_k$
is canonically %identified with 
isomorphic to $H^2(k_\xi,L)$. There is a canonical
action of the sheaf $\Cal H^2(Z)$ (of abelian groups) on the sheaf
$\Cal H^2(L)$, which is (stalkwise) free and transitive. \medbreak}

\Proof The identification of the stalks follows from Lemma (2.3).
Also it is clear from (1.24) that one has a canonical action of
$\Cal H^2(Z)$ on $\Cal H^2(L)$, which is stalkwise free and
transitive (wherever the stalks of the second sheaf are non-empty).
Since $\Cal H^2(Z)$ is known to be a locally constant sheaf with
finite stalks ([Sch], Thm.\ 2.13a), the same will follow for 
$\Cal H^2(L)$ once it is shown that the subset $\{\xi\colon$ 
$H^2(k_\xi,L)\ne\emptyset\}$ is clopen in $\Omega_k$. But this is the
set of $\xi$ where $\obs(L)_\xi=0$ in $H^3(k_\xi,Z)$, and hence it
is obviously clopen. \endproof

\n{\bf(2.5) Notation.}\enspace Given $\alpha\in H^2(k,L)$,
we write $\alpha_\xi$ for the restriction of $\alpha$ to $k_\xi$,
so $\alpha_\xi\in H^2(k_\xi,L)$. \bigbreak

\n{\bf(2.6) Corollary.\enspace\it If\/ $\vcd(k)\le1$ then
the map $H^2(k,L)\to\Gamma\bigl(\Omega_k,\,\Cal H^2(L)\bigr)$ is
bijective. 
\medbreak}

\Proof If $H^2(k,L)=\emptyset$ this is true by (1.22). Otherwise
use Corollary (2.4) together with bijectivity of $H^2(k,Z)\to\Gamma
\bigl(\Omega_k,\,\Cal H^2(Z)\bigr)$ ([Sch], Thm.\ 3.1a). \endproof

\n{\bf(2.7) Definition.}\enspace The sheaf (of sets) $\Cal N^2(L)$ 
on $\Omega_k$ is defined by
$$U\quad\longmapsto\quad\dirlim{(E,s)\in J_U}N^2(E,L),$$
for $U\subset\Omega_k$ clopen. It is clear that this is a subsheaf
of $\Cal H^2(L)$. \bigbreak

The proof that this subsheaf $\Cal N^2(L)$ of $\Cal H^2 (L)$ has
the right stalks, and that it is again locally constant requires
slightly more work than for $\Cal H^2(L)$. We first prove \bigbreak

\n{\bf(2.8) Lemma.\enspace\it Let $\alpha\in H^2(k,L)$ be an
element whose restriction to $k_\xi$ is neutral, where $\xi\in
\Omega_k$ is fixed. Then there is a finite subextension $k\subset
K\subset k_\xi$ such that already the restriction of $\alpha$ to
$K$ is neutral. \medbreak}

\Proof Let $(f,g)$ be a cocycle which represents $\alpha$. We can
assume that $(f,g)$ is normalized, i.e.\ that $f_1={\text id}$ and 
$g_{1,s}=g_{s,1}=1$ for $s\in\Gamma$. Fix a real closure $k_\xi
\subset k_s$ with respect to $\xi$, and let $t$ be the
corresponding involution in $\Gamma$. That $\alpha$ is neutral at
$\xi$ means (see $(6)$) that there is $h\in\ol G(k_s)$ such that
$$h\cdot f_t(h)\cdot g_{t,t}\>=\>1.$$
Choose an open normal subgroup $U$ of $\Gamma$ with $t\notin U$,
and put $V:=U\cup tU$. If one takes $U$ sufficiently small then the
following properties hold:
\item{a)} If $x,y\in V$ then $g_{x,y}=1$ if $x$ or $y$ is in $U$,
and $g_{x,y}=g_{t,t}$ otherwise;
\item{b)} $f_x(h)=h$ and $f_{xt}(h)=f_t(h)$ for every $x\in U$.
\par\n
b) follows since $f$ is (weakly) continuous. Now the restriction of
the cocycle $(f,g)$ to the open subgroup $V$ of $\Gamma$ is
neutralized by the continuous map $V\to\bar G(k_s)$ which is $1$ on
$U$ and $h$ on $tU$. \endproof

\n{\bf(2.9) Proposition.\enspace\it The subsheaf $\Cal N^2(L)$ of 
$\Cal H^2(L)$ is locally constant. Its stalk at $\xi$ is
the subset $N^2(k_\xi,L)$ of $H^2(k_\xi,L)$. \medbreak}

\Proof It is clear that the stalk of $\Cal N^2(L)$ at $\xi$ is
contained in $N^2(k_\xi,L)$. The other inclusion follows directly
from Lemma (2.8). A subsheaf of the locally constant sheaf $\Cal H
^2(L)$ is locally constant if and only if it is (not only open but
also) closed in $\Cal H^2(L)$, both sheaves being regarded as
espaces \'etal\'es. This follows immediately from [Sch], Lemma C.3,
which asserts that a sheaf on a boolean space is locally constant
with finite stalks if and only if its espace \'etal\'e is compact.
Therefore we have to show: Given a section $c$ of $\Cal H^2(L)$
over a clopen subset $U\subset\Omega_k$, the set of $\xi\in U$ with
$c_\xi\in N^2(k_\xi,L)$ is closed. For this in turn it suffices to
show for every finite extension $E$ of $k$ and every $\alpha\in H^2
(E,L)$ that the set $\{\eta\in\Omega_E\colon$ $\alpha_\eta$ is
neutral$\}$ is closed. But this follows again from Lemma (2.8).
\endproof
\bigskip
\subhead 3. A formally real analogue of Grothendieck's theorem\endsubhead 

\n Grothendieck's theorem, by which we mean Theorem 3.5 in [Sp],
asserts: If $k$ is a perfect field with $\cd(k)\le1$, and $L=(\ol
G,\kappa)$ is a $k$-kernel, then $H^2(k,L)$ consists of precisely
one element, and this element is neutral. In particular, any
$k$-kernel is trivial.

The following theorem specializes to Grothendieck's theorem if the
field $k$ has no orderings. It can be regarded as its ``formally
real'' analogue. 
\bigbreak

\n{\bf(3.1) Theorem.\enspace\it Let $k$ be a perfect field
with $\vcd(k)\le1$, and let $L=(\ol G,\kappa)$ be a $k$-kernel.
\item{\rm a)} The sheaf $\Cal H^2(L)$ on\/ $\Omega_k$ is locally
constant. Its stalk at $\xi$ is the finite set $H^2(k_\xi,L)$.
\item{\rm b)} The subsheaf $\Cal N^2(L)$ of $\Cal H^2(L)$ is again
locally constant. Its stalk at $\xi$ is the subset $N^2(k_\xi,L)$
of $H^2(k_\xi,L)$.
\item{\rm c)} The natural map $H^2(k,L)\to\Gamma\bigl(\Omega_k,\,
\Cal H^2(L)\bigr)$ is bijective, and its restriction to $N^2(k,L)$
is a bijection onto $\Gamma\bigl(\Omega_k,\,\Cal N^2(L)\bigr)$.
\bigbreak}

\n{\bf(3.2) Corollary.\enspace\it Let $\alpha$, $\beta\in H^2(k,L)$.
\item{\rm a)} $\alpha=\beta$ $\Leftrightarrow$ $\alpha_\xi=\beta_
\xi$ for all $\xi$ in a dense subset of\/ $\Omega_k$.
\item{\rm b)} $\alpha$ is neutral $\Leftrightarrow$ $\alpha_\xi$ is
neutral for all $\xi$ in a dense subset of\/ $\Omega_k$.
\item{\rm c)} $H^2(k,L)$ contains a neutral element
$\Leftrightarrow$ $H^2(k_\xi,L)$ contains a neutral element for all
$\xi$ in a dense subset of\/ $\Omega_k$. \medbreak}

\Proof a) and b) follow from (3.1c), using the characterization of
the stalks of the sheaves $\Cal H^2(L)$ and $\Cal N^2(L)$. c)
follows from (3.1c) together with the fact that the sheaf 
$\Cal N^2(L)$ is locally constant. \endproof

\n
Our proof does not assume Grothendieck's theorem, i.e.\ the case of
(3.1) where $k$ has no orderings; rather this case will be covered
by our proof as well. 

For the proof of Theorem (3.1) we will need the following two
lemmas. The first is the general principle underlying the main
induction step in the proof. \bigbreak

\n
\n{\bf(3.3) Lemma.\enspace\it Let $k$ be a perfect field and
let $\Cal P$ be a property of algebraic groups over $k$. (Recall
that all algebraic groups are assumed to be smooth of finite type.) 
Suppose that $\Cal P(G)$ holds under each of the following conditions 
(i) to~(iv):
\item{\rm (i)} $G$ is finite;
\item{\rm (ii)} $G$ is commutative;
\item{\rm (iii)} $G$ is connected linear semisimple with trivial
center;
\item{\rm (iv)} there is an algebraic $k$-subgroup $N$ of $G$ which
is invariant under all semilinear automorphisms of $G$ and such that 
$\Cal P(N)$ and $\Cal P(G/N)$ hold.
\par\n
Then $\Cal P(G)$ holds for every $G$. \medbreak}

\Proof Assume that the lemma is false. From the descending chain
condition for closed algebraic subschemes, there exists a minimal
element among all counterexamples. We may therefore assume that
there is a counterexample $G$ such that $\Cal P(H)$ holds for every
proper algebraic $k$-subgroup $H$ of $G$. Applying (iv) with
$N=G^0$, and using [DG], II, \S 5, 2.1(ii) (``$G$ is smooth if and 
only if $G^0$ is smooth''), (i) implies that $G$ is connected. 
Applying (iv) with $N=$ largest connected linear subgroup $L$ of $G$,
and using [R], Thm.~16, (ii) implies that $G$ is linear. Note that 
in [R] all algebraic groups are reduced over $\ol k$, hence they are
smooth over their field of definition (see [DG], II, \S 5, 2.1(v)); 
the ``$k$-closed'' subgroup $L$ of [R], Thm. 16, is defined over $k$
since $k$ is perfect, hence it is smooth over $k$. Inductively one 
also sees that $\Cal P(H)$ holds whenever $H$ is solvable. Applying 
(iv) with $N=$ radical of $G$ one sees that $S:=G/N$ is also a 
counterexample. Using (iv) with $N=$ center of $S$ one sees from (i) 
(or (ii)) and (iii) that $\Cal P(S)$ holds, thereby giving a contradiction. 
\endproof

If $H$ is any (abstract) group, let $\Inv(H)=\{h\in H\colon$
$h^2=1$, $h\ne1\}$ denote the set of involutions in $H$. The
following lemma is used in the proof of Theorem (3.1) to verify
step (iv) in Lemma (3.3), the property $\Cal P$ being the
local-global principle for neutral elements. \bigbreak

\n{\bf(3.4) Lemma.\enspace\it Let $k$ be a perfect field
with\/ $\vcd(k)\le1$, let $\ol H$ be an algebraic group over 
$\ol k=k_s$ and let
$$1\to\ol H(\ol k)\to F\rightmap\pi\Gamma\to1\eqno(10)$$
be an extension of the type considered in\/ (1.18), where\/
$\Gamma:=\Gamma_k$.
\item{\rm a)} Suppose\/ $(10)$ splits locally, i.e.\
$\Inv(F)\to\Inv (\Gamma)$ is surjective. Then there exists a
continuous map $\tau \colon\Inv(\Gamma)\to\Inv(F)$ with
$\pi\comp\tau={\text id}$ and such that for every $t\in\Inv(\Gamma)$
and $x\in F$, the elements $x\tau (t)x^{-1}$ and
$\tau\bigl(\pi(x)\cdot t\cdot\pi(x)^{-1}\bigr)$ are conjugate under
$\ol H(\ol k)$. We call such $\tau$ an\/ {\text Inv-section} of\/
$(10)$.
\item{\rm b)} If\/ $(10)$ splits, then for any Inv-section $\tau$
of\/ $(10)$ there is a splitting $\sigma\colon\Gamma\to F$ of\/
$(10)$ such that for every $t\in\Inv(\Gamma)$, $\sigma(t)$ is
conjugate to $\tau(t)$ under $\ol H(\ol k)$. \medbreak}

\Proof a) We repeatedly use the following obvious fact: If $p\colon
X\to T$ is a surjective map of topological spaces which is a local
homeomorphism, and if $T$ is a boolean space, then $p$ has a
(continuous) section. By this device there exist continuous
sections $\sigma\colon\Inv(\Gamma)\to\Inv(F)$ and $\varrho\colon
\Gamma\to F$ of $\pi$. The canonical map $\Inv(\Gamma)\to\Inv
(\Gamma)/{\text conj.}=\Omega_k$ is known to have a continuous
section ([H], Lemma 5.3(a)). Hence there is a closed subset $Z$ of
$\Inv(\Gamma)$ which is a system of representatives of conjugacy
classes in $\Inv(\Gamma)$. Let $\alpha\colon\>\Gamma\times\Inv
(\Gamma)\to\Inv(\Gamma)$ be the conjugation action of $\Gamma$ on
$\Inv(\Gamma)$. Let $\iota$ be the involution on $\Gamma\times Z$
defined by $\iota(s,t)=(st,t)$. Then $\alpha$ induces a
homeomorphism $(\Gamma\times Z)/\iota\isomap\Inv(\Gamma)$. Since
$\iota$ has no fixed points, there is a closed subset $Y$ of
$\Gamma\times Z$ such that $\alpha|Y$ is a homeomorphism $Y\isomap
\Inv(\Gamma)$. Define $\tau$ by $\tau(sts^{-1}):=\varrho(s)\sigma
(t)\varrho(s)^{-1}$ for $(s,t)\in Y$.

b) Fix a splitting $\varrho$ of $(10)$, and write ${}^xa:=\varrho
(x)\,a\,\varrho(x)^{-1}$ for $x\in\Gamma$, $a\in\ol H(\ol k)$. Let
the map $c\colon\Inv(\Gamma)\to\ol H(\ol k)$ be defined by $\tau(t)
=c(t)\varrho(t)$ for $t\in\Inv(\Gamma)$. The map $c$ is locally
constant. Using that $\tau$ is an Inv-section one finds easily that
$c$ satisfies
$$\hbox{$c(t)\cdot{}^tc(t)\>=\>1$ for $t\in\Inv(\Gamma)$}$$
and
$$\hbox{for $t\in\Inv\Gamma$ and $x\in\Gamma$ there is $a\in\ol H
(\ol k)$ with ${}^xc(x^{-1}tx)=a^{-1}\cdot c(t)\cdot{}^ta$.}$$
In other words, $c\in Z^1_{\Omega_k}(\ol H(\ol k))$, in the
terminology of [Sch], Lemma 2.7. By Prop.\ 2.9 of {\it loc.~cit.}, 
$c$ represents a global section of the sheaf $\Cal H^1(H)$, where 
$H$ denotes the $k$-form of $\ol H$ defined by $\varrho$. By [Sch],
Thm.\ 5.1 there is a class $\gamma\in H^1(k,H)$ which maps to this
global section. Let $(h_s)_{s\in\Gamma}$ be a cocycle representing
$\gamma$, and define the splitting $\sigma$ of $(10)$ by $\sigma(s)
:=h_s\cdot\varrho(s)$. Then $\sigma$ has the desired property.
\endproof

\lab(3.5) We shall now prove Theorem (3.1). Assertions a), b)
and the first part of c) have already been proved in (2.4), (2.9)
and (2.6), respectively. Therefore it remains to show the following
local-global principle for neutral elements:
\item{$(*)$} An element $\alpha\in H^2(k,L)$ whose restriction
$\alpha_\xi$ to $k_\xi$ is neutral for every $\xi\in\Omega_k$ is
neutral.
\par\n
Our proof is based on an induction over $\ol G$, using Lemma (3.3).
Take $K=\ol k$ in the lemma, and let $\Cal P(\ol G)$ be the
property that $(*)$ holds for all $k$-kernels in $\ol G$. Property
(i) required in (3.3) is known (see [Sch], Sect. 5.12). Property (ii) is
part of [Sch], Thm.\ 3.1. To prove property (iii) we remark that if
$L=(\ol G,\kappa)$ is a $k$-kernel with $\ol G$
connected reductive, then $N^2(k, L)$ is not empty. (Here $k$ may
be any field; cf.\ [D] and also [B1], Prop.\ 3.1.) If further
$\ol G$ has trivial center, then by (1.24) $H^2(k, L)$ has
exactly one element, which is neutral.

% In Sect.~4 we will give a different proof of $(*)$ which is valid
% for all connected linear kernels $L$, if $\ch(k)=0$.

It remains to show that property (iv) from (3.3) holds as well. So
let $L=(\ol G,\kappa)$ and $\alpha\in H^2(k,L)$ be given with
$\alpha_\xi\in N^2(k_\xi,L)$ for all $\xi$, and let $\ol N$ be an
algebraic subgroup of $\ol G$ which is invariant under all
semilinear automorphisms of $\ol G$ and such that $(*)$ holds for $\ol N$
and for $\ol G/\ol N$. Let $\Gamma=\Gamma_k$, and let
$$1\to\ol G(\ol k)\to E\rightmap\pi\Gamma\to1\eqno(11)$$
be the extension corresponding to $\alpha$, cf.\ (1.19). The
hypothesis says that the map $\Inv(E)\to\Inv(\Gamma)$ induced by
$\pi$ is surjective. We have to show (1.26) that $(11)$ splits. The
subgroup $\ol N(\ol k)$ of $E$ is normal, and by hypothesis we know
that $(11)$ splits modulo $\ol N(\ol k)$, i.e.\ that
$$1\to\ol G(\ol k)\big/\ol N(\ol k)\to E\big/\ol N(\ol k)\rightmap
{\bar\pi}\Gamma\to1\eqno(12)$$
splits. Since $(11)$ splits locally, there is an Inv-section $\tau$
of $(11)$, by part a) of Lemma (3.4). Since $(12)$ splits, and
since $\bar\tau:=\tau$ mod~$\ol N(\ol k)$ is an Inv-section of
$(12)$, part b) of that lemma shows that there is a splitting
$\sigma$ of $(12)$ for which $\sigma(t)$ is conjugate to
$\bar\tau(t)$ under $\ol G\big/\ol N(\ol k)$, for every
$t\in\Inv(\Gamma)$.

Fix this $\sigma$, and let $S$ be the preimage of $\sigma(\Gamma)$
under $E\to E\big/\ol N(\ol k)$. Then $S$ is a subextension
$$1\to\ol N(\ol k)\to S\to\Gamma\to1\eqno(13)$$
of $(11)$. The extension $(13)$ corresponds to a $k$-kernel in
$\ol N$. Indeed, let $s\mapsto z_s$ be a locally constant section
of $(13)$, and let $f_s$ be the (unique) $s$-semilinear automorphism of
$\ol G$ which induces $\int(z_s)$ on $\ol G(\ol k)$. Then $s\mapsto
f_s$ is weakly continuous, hence continuous by (1.13) since
$\kappa$ is a kernel. But this implies that also $s\mapsto f_s|\ol
N$ is continuous. Therefore $s\mapsto$ image of $f_s|\ol N$ in
$\SOut(\ol G/k)$ is a $k$-kernel to which $(13)$ belongs.

Since for every $t\in\Inv(\Gamma)$, $S$ contains a $\ol G(\ol
k)$-conjugate of $\tau(t)$, the extension $(13)$ splits locally. By
the induction hypothesis, therefore, $(13)$ splits, and {\it a
fortiori\/} $(11)$ splits. \endproof

An alternative approach to the above induction process, which
applies only for connected linear groups, is introduced in Sect.~4.
Using Borovoi's construction of abelianized noncommutative $H^2$,
an independent proof is given there for the local-global principle
for neutral elements, if the kernel is connected and linear. Using
this result one may in the above proof directly proceed from
general $\ol G$ to its connected component, and then to the maximal
linear connected subgroup of $\ol G$. \bigbreak

\n{\bf(3.6) Remark.}\enspace
From Theorem (3.1) it is clear that it can happen (over suitable
$k$ with $\vcd(k)=1$) that a $k$-kernel $L$ is trivial over {\it
some\/} real closure $k_\xi$, but not over another (and {\it a
fortiori\/} not over $k$), even if $H^2(k,L)\ne\emptyset$. As an
illustration we give the following explicit example: Let
$$1\to G\to E\to\Bbb Z/2\to1\eqno(14)$$
be an extension of finite groups which does not split modulo the
center of $G$. (Example: Take $E$ to be group of all
transformations $aX+b$ over the field with five elements, and $G$
the subgroup consisting of those with $a=\pm1$.) Let $\tau$ be the
involution in $\Out(G)$ defined by $(14)$. Then $\tau$ does not
lift to an involution in $\Aut(G)$.

Let $k$ be a field with $\vcd(k)=1$, and let $K/k$ be a quadratic
extension such that $K$ is formally real and there is an ordering
of $k$ which does not extend to $K$. (Example: $k=\Bbb R(t)$,
$K=k(\sqrt t)$.) Let $\chi\colon\Gamma_k\to\Bbb Z/2$ be the
homomorphism with kernel $\Gamma_K$, and let
$$1\to G\to\tilde E\to\Gamma_k\to1\eqno(15)$$
be the extension obtained by pulling back $(14)$ via $\chi$. This
extension defines a $k$-kernel $L=(G,\kappa)$ ($G$ is considered as
a constant finite group scheme here) and an element in $H^2(k,L)$.
For every ordering $\xi$ of $k$ which extends to $K$, $(15)$ is
trivial over $k_\xi$, and hence $N^2(k_\xi,L)\ne\emptyset$. On the
other hand, if $\xi$ does not extend to $K$ then $N^2(k_\xi,L)$ is
empty since $\tau$ does not lift to an involution in $\Aut(G)$.
\bigskip
\subhead 4. Hypercohomological approach\endsubhead
%\beginsect4.\ Hypercohomological approach

\n The aim of this section is to introduce another technique which 
provides an alternative proof of the induction step in Theorem
(3.1). We consider only connected linear algebraic groups in this
section. The main tool here is the abelian Galois hypercohomology
group $H^2$ of a complex of length two, which was used by Borovoi 
[B1] in the context of number fields.

\lab(4.1) Let $L=(\ol G,\kappa)$ be a connected reductive
$k$-kernel (i.e.\ $\ol G$ is a connected linear reductive algebraic
group over $k_s$). Let $\olgss=[\ol G,\ol G]$ be its derived
group; it is semisimple. We denote by $\olgsc$ the simply connected
cover of $\olgss$ and by $\ol{\rho}$ the composite map
$$\olgsc\to\olgss\to\ol G.$$
Let $\ol Z$ (resp. $\olzss$, $\olzsc$) be the center of $\ol G$
(resp. $\olgss$, $\olgsc$). Observe that $\kappa$ defines $k$-forms
$Z$,~$\zss$ and $\zsc$ of $\ol Z$,~$\olzss$ and $\olzsc$ since
these groups are abelian. Further, if $L=(G\times_k k_s,\kappa_G)$
for some $k$-group $G$ (cf. (1.15)), then the $k$-groups
$Z$,~$\zss$ and $\zsc$ are the respective centers of $G$,~$\gss$
and $\gsc$ (cf. (1.16)). The restricted homomorphism $\rho:\zsc\to
Z$ is then defined over $k$. This induces a short complex of
discrete $\Gamma_k$-modules
$$1\to\zsc(k_s)\buildrel{\rho}\over{\to}Z(k_s)\to1\eqno(16)$$
placed in degrees $-1$ and $0$. Borovoi [B1], Sect.~5, defines
the abelianized cohomology groups $H^i_{\ab}(k,L)$, $i\ge0$ by
$$H^i_{\ab}(k,L):={\Bbb H}^i(k,\,\zsc\to Z)$$
where the right hand side denotes the $\Gamma_k$-hypercohomology of
the complex $(16)$. Borovoi also constructs an abelianization map
[B1], Sect. 5.3,
$${\ab}^2:H^2(k,L)\to H^2_{\ab}(k,L).$$
These constructions are valid for a field $k$ of any
characteristic, although [B1] assumes that $\ch(k)=0$.

\lab(4.2) To extend the definition of the abelianization map to the
case of any connected linear $k$-kernel $L=(\ol G,\kappa)$ (i.e.\ 
$\ol G $ is a connected linear algebraic group over $\ol k$), we 
assume that $k$ is perfect. This assumption remains in force for
the rest of this section. 

Let $\ol G^{\red}$ denote the connected
reductive group which is the quotient of $\ol G$ by its unipotent
radical and let $L^{\red}=(\ol G^{\red},\kappa)$ be the induced
kernel. There is a natural map [B1], Sect. 4, $r:H^2(k,L)\to H^2
(k,L^{\red})$. Borovoi proved [B1], Prop. 4.1, that an element
$\eta\in H^2(k,L)$ is neutral if and only if $r(\eta)$ is neutral.
Note that [B1], Prop. 4.1, holds for any perfect field $k$. Indeed,
Lemma 4.3 of [B1] holds for all such fields, since [DG], IV, \S2,
Cor. 3.9, asserts that a (connected smooth) unipotent $k$-group
has a central composition series with quotients $\Bbb G_a$. Setting
$H^2_{\ab}(k,L)=H^2_{\ab}(k,L^{\red})$, the map $\ab^2$ is defined
as the composite
$$H^2(k,L)\buildrel{r}\over{\to}H^2(k,L^{\red})\buildrel{\ab^2}
\over{\to}H^2_{\ab}(k,L^{\red}).$$
The main advantage of the abelianization map is that it helps 
detect neutral elements. More precisely, we have \bigbreak

\n{\bf(4.3) Proposition.\enspace\it Let $k$ be a perfect
field with\/ $\vcd(k)=1$ and let $L=(\ol G,\kappa)$ be a connected
linear $k$-kernel. An element $\eta\in H^2(k,L)$ is neutral if and
only if\/ $\ab^2(\eta)$ is~$0$. \medbreak}

\Proof The proof is the same as in [B1], Sect. 5.8, once we replace [B1],
Lemma 5.7, by \medskip

\n{\bf(4.4) Lemma.\enspace\it Let $G$ be a semisimple simply
connected linear algebraic group over a perfect field $k$ with
$\vcd(k)\le1$ and put $G^{\ad}=G/Z$, where $Z$ is the center of
$G$. Then the connecting map $\delta:H^1(k,G^{\ad})\to H^2(k,Z)$ is
surjective. \medbreak}

\Proof (See [Sch], Cor. 5.4, for another proof). Since $\cd(k(i))\le1$,
Steinberg's theorem ([S2], III, \S 2.3) implies that $H^1(k(i),G)$ is
trivial and $G$ is quasi-split over $k(i)$. Therefore $G$ has a
Borel subgroup $B$ defined over $k(i)$. Further, $B$ can be chosen
such that $B\cap\sigma B$ (where $\sigma$ is the involution of
$k(i)$ over $k$) is a torus $T$ (necessarily maximal in $G$ and
defined over $k$; cf.\ [Sch], Prop. 4.9). The torus $T$ has the property
that $H^2(k_{\xi},T)=0$ (cf.\ [Sch], Prop. 1.6 and proof of Corollary
1.7) for every ordering $\xi\in\Omega_k$. Since the map $H^2(k,T)
\to\prod_\xi H^2(k_\xi,T)$ is injective [Sch], Thm. 3.1, the group $H^2
(k,T)$ is itself trivial. Using this in the long exact cohomology
sequence associated to the short exact sequence
$$1\to Z\to T\to T^{\ad}\to1$$
which defines $T^{\ad}$, we see that the map $H^1(k,T^{\ad})\to H^2
(k,Z)$ is surjective. Hence in the commutative diagram
$$\matrix H^1(k,T^{\ad})&\twoheadrightarrow&H^2(k,Z)\\
\downarrow&&\Arrowvert\\ H^1(k,G^{\ad})&\to&H^2(k,Z),\endmatrix$$
the surjectivity of the top arrow implies the surjectivity of the
bottom arrow. \endproof

The following proposition proves $(*)$ of Sect.~3 (cf.\ (3.5)) for
connected linear algebraic groups. \medskip

\n{\bf(4.5) Proposition.\enspace\it Let $L=(\ol G,\kappa)$
be a connected linear $k$-kernel, where $k$ is a perfect field with
$\vcd(k)\le1$. Then $\eta\in H^2(k,L)$ is neutral if and only if
its localizations $\eta_\xi\in H^2(k_\xi,L)$ are neutral for all
orderings $\xi$ in a dense subset of $\Omega_k$. \medbreak}

\Proof We may assume without loss of generality (cf.\ (4.2)) that
$\ol G$ is connected and reductive. Hence $H^2(k,L)$ has a
neutral element (cf.\ [B1], Prop.\ 3.1). Denote by $G$ the
corresponding $k$-form of $\ol G$. Let $T$ be a maximal $k$-torus
of $G$, $p:G^{\sc}\to G^{\ss}$ the simply connected cover of the
derived group $G^{\ss}$ of $G$, put $T^{\ss}=T\cap G^{\ss}$ and $T^
{\sc}=p^{-1}(T^{\ss})$.

In this set up, Borovoi [B2], Sect.~3, showed that the complexes
$(T^{\sc}\to T)$ and $(Z^{\sc}\to Z)$ are quasi-isomorphic.
Therefore we have
$$H^i_{\ab}(k,L)={\Bbb H}^i(k,T^{\sc}\to T)={\Bbb H}^i(k,Z^{\sc}
\to Z).$$

We consider a maximal torus $T=B\cap\sigma B$ as in the proof of
Lemma (4.4). Associated to the short exact sequence
$$1\to(1\to T)\to(T^{\sc}\to T)\to(T^{\sc}\to1)\to1$$
of short complexes placed in degrees $-1$ and $0$, there is a
commuting diagram of long exact hypercohomology sequences as
follows:
$$\matrix H^2(k,T^{\sc})&\to&H^2(k,T)&\to&H^2_{\ab}(k,L)&\to
&H^3(k,T^{\sc})\\ \downarrow&&\downarrow&&\downarrow&&\downarrow\\
\prod_\xi H^2(k_\xi,T^{\sc})&\to&\prod_\xi H^2(k_\xi,T)&\to&
\prod_\xi H^2_{\ab}(k_\xi,L)&\to&\prod_\xi H^3(k_\xi,T^{\sc}).\endmatrix$$
As noted in the proof of Lemma (4.4), the terms on the left are
zero, hence the second horizontal arrows are injections. By [Sch],
Cor. 3.2, the second and fourth vertical arrows are injective, even 
when $\xi$ ranges only over a dense subset of $\Omega_k$. Hence the
third vertical map is injective. Our proposition follows from the
commutative diagram
$$\matrix H^2(k,L)&\buildrel{\ab^2}\over{\to}&H^2_{\ab}(k,L)\\
\downarrow&&\downarrow\\
\prod_\xi H^2(k_\xi,L)&\buildrel{\ab^2_\xi}\over{\to}
&\prod_\xi H^2_{\ab}(k_\xi,L),\endmatrix$$
on using Proposition (4.3) and the injectivity of the vertical
arrow on the right. \endproof

\lab(4.6) Proposition (4.5) establishes $(*)$ of (3.5), and
consequently (iii) of Lemma (3.3), for a connected linear algebraic
group, which is not necessarily semisimple or of trivial center.
This eliminates a few steps in the proof of
% Lemma (3.3), and consequently in the proof of
Theorem (3.1). \bigskip
\bigskip
\subhead 5. Applications\endsubhead 

\n As an application of Grothendieck's theorem (more precisely, of the
local-global principle for neutral elements) we get a new proof of
Theorem 6.1 of [Sch]. We can also answer now the question raised in
{\it loc.cit.} after 6.5.

\lab(5.1) Let $k$ be a field, let $G$ be an algebraic group over
$k$ and let $X$ be a homogeneous (right) $G$-space, defined over
$k$. From $X$ one constructs canonically a $k$-kernel $L_X$ and an
element $\alpha_X\in H^2(k,L_X)$. The class $\alpha_X$ is neutral
if and only if there is a principal homogeneous $G$-space $T$ which
dominates $X$ (over $k$).

Let us recall this construction ([Sp], Sect. 1.20; [B1], Sect. 7.7). 
Choose a point $x_0\in X(k_s)$, and let $\ol H$ be the stabilizer of 
$x_0$. This is an algebraic subgroup of $\ol G:=G\times_kk_s$. There 
is a locally constant map $\Gamma=\Gamma_k\to G(k_s)$, $s\mapsto a_s$,
such that $s(x_0)=x_0 \cdot a_s$ for every $s\in\Gamma$. Then
$\int(a_s)\comp s$ is an $s$-semilinear automorphism of $\ol G$ which
leaves $\ol H$ invariant. Let $f_s$ be its restriction to $\ol H$;
the map $s\mapsto f_s$ from $\Gamma$ to $\SAut(\ol H/k)$ is
continuous (1.10), and $s\mapsto f_s$ mod~$\Int(\ol H)$ is a
$k$-kernel in $\ol H$. Denote it by $\kappa_X$, and write $L_X=(\ol
H,\kappa_X)$. Up to a canonical isomorphism, this kernel does not
depend on the choices. The pair $(f,g)$ with $f$ as above and
$g_{s,t}:=a_ss(a_t)a_{st}^{-1}$ is a 2-cocycle for $L_X$, and
$\alpha_X$ is its class in $H^2(k,L_X)$. Again, this class does not
depend on the choices. It is neutral if and only if $X$ is
dominated by a principal homogeneous $G$-space, as one may verify
by a direct cocycle calculation.

In terms of group extensions the class $\alpha_X$ has a more
natural appearance: Let $G(k_s)\cdot\Gamma$ be the natural
semidirect product, and let $E=E_{X,x_0}$ be its subgroup
consisting of all products $gs$ with $g\in G(k_s)$, $s\in\Gamma$
such that $s(x_0)=x_0\cdot g$. Then $E$ is a subextension of the
split extension $G(k_s)\cdot\Gamma$ as follows:
$$1\to\ol H(k_s)\to E\to\Gamma\to1.\eqno(17)$$
This extension belongs to the $k$-kernel $\kappa_X$, and the class
of $(17)$ in $H^2(k,L_X)$ is $\alpha_X$. \bigbreak

\n{\bf(5.2) Theorem.\enspace\it Let $k$ be a perfect field
with $\vcd(k)\le1$. Let $X$ be a homogeneous space under an
algebraic group $G$, both defined over $k$, and suppose that for
every $\xi$ in a dense subset of\/ $\Omega_k$, $X$ is dominated by
some principal homogeneous $G$-space over $k_\xi$. Then there
exists a principal homogeneous $G$-space over $k$ which dominates
$X$ over $k$. \medbreak}

\n This was proved in [Sch], Thm.\ 6.5, under the stronger
hypothesis that $X$ is dominated by a principal homogeneous space
over {\it every\/} $k_\xi$, and the question was raised there
(p.~341) whether the above sharpening holds. The proof is very easy
now: \medbreak

\Proof Let the kernel $L_X$ and the class $\alpha_X\in H^2(k,L_X)$
be as constructed in (5.1). By the above, the hypothesis says that
the class $\alpha_X$ becomes neutral over $k_\xi$ for all $\xi$ in
a dense subset. By Corollary (3.2), $\alpha_X$ is neutral over $k$.
Again by (5.1), this translates into the existence of a principal
homogeneous $G$-space as desired. \endproof

In [Sch], the above-mentioned Thm.\ 6.5 was derived from the
following main result: \bigbreak

\n{\bf(5.3) Theorem\enspace\rm([Sch], Thm.\ 6.1).\enspace\it
Let $k$ be a perfect field with $\vcd(k)\le1$, and let $X$ be a
homogeneous space under an algebraic group $G$, both defined over
$k$. Suppose $X(k_\xi)\ne\emptyset$ for all $\xi$ in a dense subset
of\/ $\Omega_k$. Then there exists a principal homogeneous
$G$-space $T$ over $k$ which dominates $X$ over $k$ and is trivial
over every $k_\xi$. \medbreak}

\n
We show how conversely Theorem (5.3) can be derived from Theorem
(5.2) (and thus from our main results in Sect.~3); in fact, the
weaker version of (5.2) proved in [Sch] suffices for this.

First of all we have $X(k_\xi)\ne\emptyset$ for all orderings
$\xi$, by general reasons (e.g.\ [Sch], Cor.\ 2.2). Therefore by
[Sch], Thm.\ 6.5 --- cf.\ (5.2) above --- there exists a $G$-torsor
$P$ together with a $G$-equivariant map $\alpha\colon\>P\to X$,
both defined over $k$. Let $H$ be the algebraic group over $k$
which consists of those $G$-equivariant automorphisms of $P$ which
commute with $\alpha$. For any extension $K/k$ the set $H^1(K,H)$
parametrizes the $K$-forms $\beta\colon\>Q\to X$ of $\alpha\colon\>
P\to X$. So we have a commutative diagram
$$\matrix H^1(k,H)&\to&H^1(k,G)\\ \downarrow&&\downarrow\\
\Gamma(\Omega_k,\,\Cal H^1(H))&\to&\Gamma(\Omega_k,\,\Cal H^1(G))
\endmatrix\eqno(18)$$
in which the upper horizontal arrow takes $\beta\colon\>Q\to X$ to
$Q$.

For every $\xi$ let $\gamma_\xi\in H^1(k_\xi,H)$ be the class of a
$G$-map $T_\xi\to X$ over $k_\xi$, where $T_\xi$ is the trivial
$G$-torsor over $k_\xi$. The family $(\gamma_\xi)_\xi$ is locally
constant, i.e.\ is a global section of the sheaf $\Cal H^1(H)$.
% (!)
Since the vertical maps in $(18)$ are surjective by [Sch], Thm.\
5.1, there exists a $G$-torsor $\beta\colon\>Q\to X$ over $X$ which
is trivial over each $k_\xi$. \endproof

Keep the assumption that $k$ is perfect with $\vcd(k)\le1$, and 
suppose that $G$ is connected and linear. Recall that any $G$-torsor 
which is trivial over $k_\xi$ for all $\xi$ in a dense subset of 
$\Omega_k$ is trivial over $k$ ([Sch], Cor. 4.2). Combining this with 
Theorem (5.3) immediately implies the following Hasse principle for 
general homogeneous spaces under $G$: \bigbreak

\n{\bf(5.4) Corollary\enspace\rm([Sch], Cor.\ 6.2).\enspace
\it Let $k$ be a perfect field with\/ $\vcd(k)\le1$, $G$ a
connected linear algebraic group and $X$ a homogeneous $G$-space,
both defined over $k$. If $X$ has a $k_\xi$-point for all $\xi$ in
a dense subset of\/ $\Omega_k$, then $X$ has a $k$-point.
\medbreak}

As noted after (6.2) of [Sch], the conclusion holds for any $G$ for
which the natural map $H^1(k,G)\to\prod_\xi H^1(k_\xi,G)$ has
trivial kernel.
\bigskip
\n{\bf References}\medskip
\item{[BP1]} E.~Bayer-Fluckiger, R.~Parimala, Galois cohomology of
the classical groups over fields of cohomological dimension $\le2$.
Invent.\ math.\ {\bf122}, 195-229 (1995).
\item{[BP2]} E.~Bayer-Fluckiger, R.~Parimala, Classical groups and
Hasse principles. Ann.\ Math., to appear.
\item{[BS]} A.~Borel, J.-P.~Serre, Th\'eor\`emes de finitude en
cohomologie galoisienne. Comment.\ Math.\ Helv.\ {\bf39}, 111-164
(1964).
\item{[B1]} M.~V.~Borovoi, Abelianization of the second nonabelian
Galois cohomology. Duke Math.\ J. {\bf72}, 217-239 (1993).
\item{[B2]} M.~V.~Borovoi, Abelian Galois cohomology of reductive
groups. Memoirs AMS\ {\bf626}, (1998).
\item{[Br]} L.~Breen, Tannakian categories. In: Motives.
U.~Jannsen, S.~Kleiman, J.-P.~Serre (eds.), Proc.\ Symp.\ Pure
Math.\ {\bf55}, Part~I, Providence, R.I., 1994, pp.\ 337-376.
\item{[CT]} J.-L.~Colliot-Th\'el\`ene, Groupes lin\'eaires sur les
corps de fonctions de courbes r\'eelles. J. reine angew.\ Math.\
{\bf474}, 139-167 (1996).
\item{[DM]} P.~Deligne, J.~Milne, Appendix to: Tannakian
categories. Lect.\ Notes Math.\ {\bf900}, 220-226 (1982).
\item{[DG]} M.~Demazure, P.~Gabriel, Groupes Alg\'ebriques,
North-Holland Publ.\ Co., Amsterdam, 1970.
\item{[D]} J.-C.~Douai, Cohomologie galoisienne des groupes
semi-simples d\'efinis sur les corps globaux. C.~R. Acad.\ Sc.\
Paris {\bf281}, 1077-1080 (1975).
\item{[FJ]} G.~Frey, M.~Jarden, Approximation theory and the rank
of abelian varieties over large algebraic fields. Proc.\ London
Math.\ Soc.\ (3) {\bf28}, 112-128 (1974).
\item{[G]} J.~Giraud, Cohomologie non ab\'elienne. Grundlehren der
mathematischen Wissen\-schaf\-ten {\bf179}, Springer, Berlin, 1971.
\item{[H]} D.~Haran, Closed subgroups of $G(\Bbb Q)$ with
involutions. J.~Algebra {\bf129}, 393-411 (1990).
\item{[M]} S.~Mac~Lane, Homology. Grundlehren der mathematischen
Wissenschaften {\bf114}, Springer, Berlin, 1963.
\item{[R]} M.~Rosenlicht, Some basic theorems on algebraic groups.
Am.\ J. Math.\ {\bf78}, 401-443 (1956).
\item{[Scha]} W.~Scharlau, Quadratic and Hermitian Forms.
Grundlehren der mathematischen Wissenschaften {\bf270}, Springer,
Berlin, 1985.
\item{[Sch]} C.~Scheiderer, Hasse principles and approximation
theorems for homogeneous spaces over fields of virtual
cohomological dimension one. Invent.\ math.\ {\bf125}, 307-365
(1996).
\item{[S1]} J.-P.~Serre, Groupes alg\'ebriques et corps de classes.
Hermann, Paris, 1959.
\item{[S2]} J.-P.~Serre, Cohomologie Galoisienne. Cinqui\`eme
\'edition. Lect.\ Notes Math.\ {\bf5}, Springer, Berlin, 1994.
%\item{[S3]} J.-P.~Serre, Cohomologie Galoisienne: Progr\`es et
%Probl\`emes. S\'em. Bourbaki 783, Mars 1994.
\item{[Sp]} T.~A.~Springer, Nonabelian $H^2$ in Galois cohomology.
In: Algebraic Groups and Discontinuous Subgroups, ed.\ A.~Borel,
G.~D.~Mostow, Proc.\ Symp.\ Pure Math.\ IX, Providence, R.I., 1966,
pp.\ 164-182.
\medskip
\n
Department of Mathematics, The Ohio State University, 231 W.~18th
Ave., Columbus, OH 43210-1174, USA; e-mail: {\tt
flicker\@math.ohio-state.edu}
\smallskip\n
Fakult\"at f\"ur Mathematik, Universit\"at Regensburg, 93040
Regensburg, Germany;\hfill\break
e-mail: {\tt claus.scheiderer\@mathematik.uni-regensburg.de}
\smallskip\n
School of Mathematics, Tata Institute of Fundamental Research, Homi
Bhabha Road, Colaba, Bombay $400\,005$, India; e-mail: {\tt
sujatha\@math.tifr.res.in}
\enddocument
\bye